\documentclass[11pt]{amsart}
\usepackage{amsmath,amssymb,amsfonts,amsthm,latexsym}

\font\tenrsfs=rsfs7 at 10pt
\font\sevenrsfs=rsfs7
\font\fiversfs=rsfs5
\newfam\rsfsfam\textfont\rsfsfam=\tenrsfs\scriptfont\rsfsfam=\sevenrsfs\scriptscriptfont\rsfsfam=\fiversfs

\def\curly{\fam\rsfsfam\tenrsfs} 

\newcommand{\PP}{{\curly P}}
\newcommand{\Cal}{\mathcal}

\newtheorem{thm}{Theorem}
\newtheorem{lem}{Lemma}[section]
\newtheorem{cor}[thm]{Corollary}

\newcommand{\N}{\Cal{N}}

\renewcommand{\a}{\ensuremath{\alpha}}

\renewcommand{\b}{\ensuremath{\beta}}

\newcommand{\lam}{\ensuremath{\lambda}}
\newcommand{\eps}{\ensuremath{\varepsilon}}
\renewcommand{\th}{\ensuremath{\vartheta}}

\newcommand{\bo}{\overline{\Omega}}

\newcommand{\fl}[1]{{\ensuremath{\left\lfloor #1 \right\rfloor}}}

\newcommand{\pfrac}[2]{{\left(\frac{#1}{#2}\right)}}

\def\({\left(}
\renewcommand{\)}{\right)}

\newcommand{\be}{\begin{equation}}
\newcommand{\ee}{\end{equation}}
\newcommand{\benn}{\begin{equation*}}   
\newcommand{\eenn}{\end{equation*}}
\newcommand{\bal}{\begin{align*}}
\newcommand{\ea}{\end{align*}}
\newcommand{\eal}{\ensuremath{\end{align*}}}
\newcommand{\bea}{\begin{eqnarray}}
\newcommand{\eea}{\end{eqnarray}}

\numberwithin{equation}{section}

\renewcommand{\le}{\leqslant}

\renewcommand{\ge}{\geqslant}

\renewcommand{\d}{\,{\rm d}}
\newcommand{\er}{{\rm e}}
\newcommand{\dsp}{\displaystyle}
\newcommand{\bline}{\noalign{\vskip-3mm}}
\def\HH{{\Cal H}}
\def\ft#1#2{\textstyle\frac{#1}{#2}}

\def\sct{\scriptstyle}
\def\di#1#2{\substack{\sct#1 \\ \noalign{\vskip-2mm}\\ \sct#2}}
       \countdef\temps=170
        \temps=\time
        \countdef\nminutes=171{\nminutes = \time}
        \countdef\nheure=172
        \def\heure{\begingroup                     
           \temps = \time \divide\temps by 60
           \nheure = \temps                        
           \nminutes = \time
           \multiply\temps by 60
           \advance\nminutes by -\temps            
           \ifnum\nminutes<10 \toks1 = {0}%
           \else\toks1 = {}%
           \fi
           \number\nheure h\the\toks1 \number\nminutes
        \endgroup}%

        \def\dater{\vglue-10mm\rightline{(\the\day/\the\month/\the\year)}}
        \def\dateheure{\the\day/\the\month/\the\year,\ \heure}

\begin{document}

\title{The distribution of integers with\goodbreak
at least two divisors in a short interval}

\author{Kevin Ford
     \and G\'erald Tenenbaum}
\thanks{First author supported by National Science
        Foundation grant DMS-0301083.}

\address{
Kevin Ford: Department of Mathematics, University of Illinois at
   Urbana-Champaign,\goodbreak 1409 West Green St., Urbana, IL, 61801, USA 
\goodbreak
G\'erald Tenenbaum:
Institut \'Elie Cartan, Universit\'e Henri--Poincar\'e Nancy 
1, B.P. 239, \goodbreak
54506 Vand\oe uvre-l\`es-Nancy Cedex, France}

\email{ford@math.uiuc.edu \quad
gerald.tenenbaum@iecn.u-nancy.fr}
\date{9 February 2007}

\begin{abstract}
We estimate the density of integers which have more than one divisor 
in an interval $(y,z]$ with $z\approx y + y/(\log y)^{\log 4 -1}$.
As a consequence, we determine the precise range of $z$ such that
most integers which have at least one divisor in $(y,z]$ have exactly
one such divisor.
\end{abstract}

\maketitle

%
\section{Introduction}\label{sec:intro}
%
Whereas, in usual cases, sieving by a set of primes may be fairly
well controlled, through Buchstab's identity, sieving
by a set of integers is a much more complicated task. However, some
fairly precise results are known in the case where
the set of integers is an interval. We refer to the recent
work~\cite{F} of the first author for  specific statements
and references.
\par
Define
\bal
\tau(n;y,z) &:= | \{ d|n : y<d\le z \} |, \\
H(x,y,z) &:= |\{ n\le x: \tau(n;y,z) \ge 1 \}|, \\
H_r(x,y,z) &:= |\{ n\le x: \tau(n;y,z) = r \}|, \\
H_2^*(x,y,z) &:= |\{ n\le x: \tau(n;y,z) \ge 2 \}| = \sum_{r\ge 2} H_r(x,y,z).
\end{align*}
\par
Thus, the numbers $H_r(x,y,z)$ $(r\ge 1)$ describe the local laws of
the function $\tau(n;y,z)$. When $y$ and $z$ are
close, it is expected that, if an integer has at least a divisor in
$(y,z]$, then it usually has exactly one,
in other words\par
\be
\label{HH1}
H(x,y,z)\sim H_1(x,y,z).
\ee
In this paper, we address the problem of determining the exact range
of validity of such behavior. In other words, we
search  for a necessary and sufficient condition so that
$H_2^*(x,y,z)=o\big(H(x,y,z)\big)$ as $x$ and $y$ tend to
infinity.
We show below that \eqref{HH1} holds if and only if
$$
\lfloor y \rfloor +1 \le z < y + \frac{y}{(\log y)^{\log 4-1+o(1)}}
  \qquad (y\to \infty).
$$
As with the results in \cite{F},
the ratios $H(x,y,z)/x$ and $H_r(x,y,z)/x$ are weakly dependent on $x$
when $x\ge y^2$.  We take pains to prove results which are valid
throughout the range $10 \le y\le \sqrt{x}$, since many interesting
applications require bounds for $H(x,y,z)$ and $H_r(x,y,z)$ when
$y\approx \sqrt{x}$; see e.g. \S 1 of \cite{F} and Ch. 2 of
\cite{Divisors} for some examples.
\par  As shown in
\cite{Te84}, for given
$y$, the threshold for the behavior of the function
$H(x,y,z)$  lies near the critical value
$$
z=z_0(y) := y\exp\{(\log y)^{1-\log 4} \} \approx y+y/(\log y)^{\log 4-1}.
$$
We concentrate on the case $z_0(y) \le z \le \er y$.  Define
\begin{equation}
\label{alldefs}
\begin{aligned}
z = \er^{\eta}& y, \quad \eta = (\log y)^{-\b},\quad
\b = \log 4 - 1 - \Xi/\sqrt{\log_2 y},\quad  \lam=\frac{1+\b}{\log 2},\\
Q(w) &= \int_1^{w} \log t\d t = w \log w - w + 1.
\end{aligned}
\end{equation}
Here $\log_k $ denotes the $k$th iterate of the
logarithm.\goodbreak\par\goodbreak
With the above notation, we have
$$\log (z/y)=\frac{\er^{\Xi\sqrt{\log_2y}}}{(\log y)^{\log 4-1}},\qquad \log
(z/z_0(y))=\frac{\er^{\Xi\sqrt{\log_2y}}-1}{(\log y)^{\log 4-1}},$$
so
\bea
\label{boundsXi}
&0\le \Xi\le (\log
4-1)\sqrt{\log_2y},\\
&0 \le \beta\le \log 4-1,\\
&\dsp\frac1{\log 2} \le \lambda\le 2.
\eea
    From Theorem 1 of \cite{F}, we know that, uniformly in $10 \le y\le
\sqrt{x}$,
$z_0(y) \le z\le \er y$,
\be\label{ford1}
H(x,y,z) \asymp
\frac{\b x}{(\Xi+1) (\log y)^{Q(\lam)}}\cdot
\ee
By Theorems 5 and 6 of \cite{F}, for any $c>0$ and uniformly in
$y_0(r) \le y \le x^{1/2-c}$, $z_0(y) \le z\le \er y$ for a suitable 
constant $y_0(r)$, we have
\be\label{ford2}
\begin{split}
\frac{H_1(x,y,z)}{H(x,y,z)} &\asymp_c 1, \\
\frac{\Xi+1}{\sqrt{\log_2 y}} \ll_{r,c} \frac{H_r(x,y,z)}{H(x,y,z)} &\le
1 \quad (r\ge 2). \\
\end{split}
\ee

When $0 \le \Xi \le o(\sqrt{\log_2 y})$ and $r\ge 2$, the upper and
lower bounds above for $H_r(x,y,z)$ have different orders.  We show
in this paper that the lower bound represents the correct order of
magnitude.

\begin{thm}\label{main_thm}
Uniformly in $10 \le y \le \sqrt{x}$, $z_0(y) \le z \le \er y$, we
have
$$
\frac{H_2^*(x,y,z)}{H(x,y,z)} \ll \frac{\Xi+1}{\sqrt{\log_2 y}}
$$
where $\Xi=\Xi(y,z)$ is defined as in \eqref{alldefs} and therefore 
satisfies \eqref{boundsXi}.
\end{thm}

\begin{cor}\label{cor_Hr}  Let $r\ge 2$ and $c>0$. There exists a 
constant $y_0(r,c)$ such that,
uniformly for $y_0(r,c) \le y \le
x^{1/2-c}$,  $z_0(y) \le z \le \er y$, we have
$$
\frac{H_r(x,y,z)}{H(x,y,z)} \asymp_{r,c} \frac{\Xi+1}{\sqrt{\log_2 y}}\cdot
$$
\end{cor}

Theorem \ref{main_thm} tells us that
$H_2^*(x,y,z)=o\big(H(x,y,z)\big)$ whenever $z\ge z_0(y)$
and $\Xi=o\big(\sqrt{\log_2 y}\big)$.
It is a simple matter to adapt the proofs given in \cite{HT91}
to show that this latter
relation persists in the range $\lfloor y \rfloor +1
\le z \le z_0(y)$. We thus obtain the following statement.

\begin{cor}\label{cor_HH1}
If $y\to\infty$, $y \le \sqrt{x}$, and
$\lfloor y \rfloor+1 \le z \le y + y(\log y)^{1-\log 4 + o(1)}$, we have
$$
H_1(x,y,z) \sim H(x,y,z).
$$
\end{cor}

Since we know from \eqref{ford2} that $H_2^*(x,y,z)\gg_\eps H(x,y,z)$
when $\beta\le \log 4-1-\varepsilon$ for any fixed $\varepsilon>0$ we have
therefore completely answered the question raised
at the beginning of this introduction concerning the exact validity
range for the asymptotic formula \eqref{HH1}. This
may be viewed as a complement to a theorem of Hall (see
\cite{Hall96}, ch.$\thinspace$7;
following a note mentioned by Hall in private correspondence, we
slightly modify the statement)
according to which
\be\label{Hall}
H(x,y,z)\sim F(-\Xi)\sum_{r\ge 1}rH_r(x,y,z)=F(-\Xi)\sum_{n\le x}\tau(n;y,z)
\ee
in the range $\Xi=o(\log_2y)^{1/6}$, $x>\exp\{\log z\log_2z\}$ with
$$F(\xi):={\frac{1}{\sqrt{\pi}}}\int_{-\infty}^{\xi/\log 4}\er^{-u^2}\d u. $$
It is likely that \eqref{Hall} still holds in the range
$(\log_2y)^{1/6}\ll \Xi\le o\big(\sqrt{\log_2y}\big)$.

%
\section{Auxiliary estimates}
%
In the sequel, unless otherwise indicated, constants implied by 
Landau's $O$ and  Vinogradov's $\ll$
symbols are absolute and effective. Numerical values of reasonable 
size could easily be given if needed.
\par
Let $m$ be a positive integer. We denote by $P^-(m)$ the
smallest, and by $P^+(m)$   the
largest, prime factor of $m$,  with the convention that
$P^-(1)=\infty$, $P^+(1)=1$.
We write $\omega(m)$ for the number of distinct prime factors of $m$ and
$\Omega(m)$ for the number of prime power divisors of $m$.  We further define
$$
\omega(m;t,u) = \sum_{\substack{ p^\nu \| m \\ t<
p\le u}} 1,\quad \Omega(m;t,u) = \sum_{\substack{ p^\nu \| m \\ t<
p\le u}} \nu,  \quad
\bo(m;t) =\Omega(m;2,t), \quad  \bo(m) = \Omega(m;2,m).
$$
Here and in the sequel, the letter $p$ denotes a prime number. Also, 
we let $\PP(u,v)$ denote the set of integers all of
whose prime factors are
in $(u,v]$ and write $\PP^*(u,v)$ for the set of squarefree members
of $\PP(u,v)$.
By convention, $1\in \PP^*(u,v)$.

\begin{lem}\label{sum1p}  There is an absolute constant $C>0$ so that for
$\ft32 \le u < v$, $v\ge \er^4$, $0 \le \a \le 1/\log v$, we have
$$
\sum_{\substack{m\in\PP(u,v) \\ \omega(m)=k}}
\frac{1}{m^{1-\a}} \le \frac{(\log_2 v - \log_2 u + C)^k}{k!}\cdot
$$
\end{lem}

\begin{proof}  For a prime $p \le v$, we have $p^\a \le 1 + 2\a \log p$, thus
the sum in question is
$$
\le \frac{1}{k!} \biggl( \sum_{u< p\le v} \frac{1}{p^{1-\a}} +
\frac{1}{p^{2-2\a}} + \cdots \biggr)^k \le \frac{\{\log_2 v - \log_2 u +
      O(1)\}^k}{k!}.
$$
\end{proof}
We note incidentally that a similar lower bound is available when $u$ 
and $v$ are not too close. See for instance Lemma
III.13 of
\cite{HR83}.

\begin{lem}\label{lemA}
Uniformly for $u\ge 10$, $0 \le k \le 2.9 \log_2 u$, and $0\le \a \le
1/(100\log u)$,  we have
$$
\sum_{\substack{P^+(m) \le u \\ \bo(m)=k}} \frac{1}{m^{1-\a}}
\ll \frac{(\log_2 u)^k}{k!}.
$$
\end{lem}

\begin{proof}
We follow the proof of Theorem 08 of \cite{Divisors}.  Let $w$ be a complex
number with $|w| \le \ft{29}{10}$.  If $p$ is prime and $3\le p \le u$, then
$|w/p^{1-\a}| \le \ft{99}{100}$ and $p^\a\le 1 + 2\a \log p$.  Thus,
$$
S(w) := \sum_{P^+(m)\le u} \frac{w^{\bo(m)}}{m^{1-\a}} = \( 1 -
\frac{1}{2^{1-\a}} \)^{-1} \prod_{3\le p\le u} \( 1 - \frac{w}{p^{1-\a}}
\)^{-1} \ll \er^{(\Re w) \log_2 u}.
$$
Put $r:=k/\log_2 u$.  By Cauchy's formula and
Stirling's formula,
\begin{align*}
\sum_{\substack{P^+(m) \le u \\ \bo(m)=k}} \frac{1}{m^{1-\a}} &= \frac{1}{2\pi
      r^k} \int_{-\pi}^{\pi} \er^{-ik\th} S(r\er^{i\th})\, d\th \ll 
\frac{(\log_2 u)^k}{k^k} \int_{-\pi}^{\pi}
\er^{k\cos \th}\, d\th \ll
\frac{(\log_2 u)^k}{k!}.
\end{align*}
\end{proof}

\begin{lem}\label{lem23}  Suppose $z$ is large, $0\le a+b\le \ft52\log_2 z$ and
$$
\exp \{ (\log x)^{9/10} \} \le w \le z\le x, \qquad x z^{-1/(10\log_2 z)} \le
Y \le x.
$$
The number of integers $n$ with $x-Y < n\le x$, $\bo(n;w)=a$
and $\omega(n;w,z)=\Omega(n;w,z)=b$, is
$$
\ll \frac{Y}{\log z} \frac{\{\log_2 w\}^a}{a!} \frac{(b+1) \{\log_2
      z - \log_2 w + C\}^b}{b!},
$$
where $C$ is a positive absolute constant.
\end{lem}

\begin{proof}  There are $\ll x^{9/10}$ integers with $n\le x^{9/10}$ or
$2^j|n$   with $2^j \ge x^{1/10}$.  For other $n$, write $n=rst$, where
$P^+(r) \le w$, $s\in \PP^*(w,z)$ and $P^-(t)>z$.  Here
$\bo(r)=a$ and $\omega(s)=b$.  We have either $t=1$ or $t>z$. In the
latter case
$x/rs>z$, whence $Y/rs>\sqrt{z}$. We may therefore apply a standard
sieve estimate to bound, for given $r$ and $s$,
the number of
$t$ by
$$
\ll \frac{Y}{rs \log z}\cdot
$$
By Lemmas \ref{sum1p} and \ref{lemA},
$$
\sum_{r,s} \frac{1}{rs} \ll \frac{(\log_2 w)^a (\log_2 z - \log_2 w
+ C)^b}{a! b!}.
$$
\par
If $t=1$, then we may assume $a+b\ge 1$.  Set $p=P^+(n)$.
If $b\ge 1$, then $p|s$ and we put $r_1:=r$ and $s_1:=s/p$.
Otherwise, let $r_1:=r/p$ and $s_1:=s=1$.  Let $A:=\bo(r_1)$ and
$B:=\omega(s_1)$, so that $A+B= a+b-1$ in all
circumstances. We have
$$
p\ge x^{1/2\bo(n)} \ge x^{1/5\log_2 z} \ge (x/Y)^2.
$$
Define the non-negative integer $h$ by $z^{\er^{-h-1}} < p \le
z^{\er^{-h}}$.  By the
Brun-Titchmarsh theorem, we see that, for each given $h$, $r_1$ and $s_1$,
the number of $p$  is
$\ll Y \er^h/(r_1 s_1\log z)$.  Set $\a:=0$ if $h=0$ and
$\a:=\er^h/(100\log z) $ otherwise.  For $h\ge 1$, we have
$r_1s_1>x^{3/4}z^{-1/\er}>\sqrt{z}$. Therefore, for $h\ge 0$,
$$
\frac{1}{r_1s_1} \le \frac{z^{-\a/2}}{(r_1s_1)^{1-\a}} \ll
\frac{\er^{-\er^{h}/200}}{(r_1s_1)^{1-\a}}.
$$
Now, Lemmas \ref{sum1p} and \ref{lemA} imply that
\begin{align*}
\sum_{r_1, s_1} \frac{1}{(r_1s_1)^{1-\a}} &\ll
\frac{(\log_2 w)^A (\log_2 z-\log_2 w+C)^B}{A! B!} \\
&\ll (b+1) \frac{(\log_2 w)^a (\log_2 z-\log_2 w+C)^b}{a!b!},
\end{align*}
where we used the fact that $a\ll\log_2w$. Summing over all $h$, we
derive that the number of those integers $n>x^{9/10}$
satisfying the conditions of the statement is
$$
\ll \frac{Y}{\log z} (b+1) \frac{(\log_2 w)^a (\log_2 z-\log_2
      w+C)^b}{a!b!}.
$$
Since $a! b! \le (3\log_2 z)^{3\log_2 z}$, this last expression is
$>x^{9/10}$. This completes the proof.
\end{proof}


Our final lemma is a special case of a theorem of Shiu (Theorem 03 of
\cite{Divisors}).

\begin{lem}\label{Shiu}
Let $f$ be a multiplicative function such that $0\le f(n) \le 1$ for all $n$.
Then, for all $x,$~$Y$ with $1< \sqrt{x} \le Y \le x$, we have
$$
\sum_{x-Y < n \le x} f(n) \ll \frac{Y}{\log x} \exp \biggl\{
\sum_{p\le x} \frac{f(p)}{p} \biggr\}.
$$
\end{lem}

%
\section{Decomposition and outline of the proof}
%

Throughout, $\eps$ will denote a very small positive constant.
Note that Theorem \ref{main_thm}
holds trivially for
$\beta \le \log 4 - 1 -
\eps$ since we then have $1\ll\Xi/\log_2y$ and of course 
$H_2^*(x,y,z)\leqslant H(x,y,z)$. We may henceforth assume that
\be\label{zrange}
\log 4 - 1 - \eps \le \beta \le \log 4 -1.
\ee
Let $$K: = \lfloor \lam \log_2 z \rfloor,$$
so that $(2-\ft32\eps)\log_2z\le K\le 2\log_2z$.
In light of \eqref{ford1}, Theorem \ref{main_thm} reduces to
\be\label{goal}
H_2^*(x,y,z) \ll \frac{x}{(\log y)^{Q(\lam)} \sqrt{\log_2 y}}\cdot
\ee
At this stage, we notice for further reference that, by Stirling's
formula, for $k\le K$ we have
\be\label{factorial}
\frac{\eta (2\log_2 z)^k}{k! (\log z)^2} \le \frac{\eta(2\log_2 z)^{K}}{K!
      (\log z)^2}
\asymp \frac{1}{(\log y)^{Q(\lam)} \sqrt{\log_2 y}}\cdot
\ee

Let $\Cal H$ denote the set of integers $n\le x$ with $\tau(n;y,z)\ge 2$.
We count separately the integers $n\in \Cal H$
lying in 6 classes.  In these definitions, we write $k=\bo(n;z)$,
$b=K-k$ and for brevity we put $z_h = z^{\er^{-h}}$.
Let
$$K_0:=(2-3\eps)\log_2z
$$ and define
\begin{align*}
\N_0 &:= \{ n\in \Cal H : n \le x/\log z \text{ or }\exists d> \log
z: d^2|n \}, \\
\N_1 &:= \{ n\in \Cal H\smallsetminus\N_0:  k \notin (K_0,K]\}, \\
\N_2 &:=\bigcup_{1\le h\le 5\eps \log_2 z}\N_{2,h},\\
\hbox{with }\N_{2,h}&:= \Big\{n\in \Cal
H\smallsetminus(\N_0\cup\N_1):\bo(n;z_h,z) \le
     \ft{19}{10}h-\ft1{100}b \Big\}.
\end{align*}
For integers $n\in \N_2$, we will only use the fact that $\tau(n;y,z)\ge 1$.
Integers in other classes do not have too many small prime factors
and it is sufficient to count pairs of divisors
$d_1,d_2$ of $n$ in $(y,z]$.  For each such pair, write
$v=(d_1,d_2)$, $d_1=vf_1$, $d_2=vf_2$, $n=f_1f_2v u$
and assume $f_1<f_2$.  Let
\be
\label{notnpf}
F_1 = \bo(f_1), \quad F_2 = \bo(f_2), \quad V= \bo(v), \quad U=\bo(u,z),
\ee
and
\be
\label{Z}
Z: = \exp \{ (\log z)^{1-4\eps} \}.
\ee
For further reference, we note that if $n\not\in\N_0$ and $h\le 5\eps\log_2
z$, then
$$
\bo(n;z_h,z)=\omega(n;z_h,z).
$$
Now we define $\HH^*:=\Cal H \smallsetminus (\N_0 \cup \N_1 \cup \N_2)$ and
\begin{align*}
\N_3 &:=  \{ n\in\HH^*  :
\min(u,f_2) \le Z  \}, \\
\N_4 &:=  \{ n\in \HH^* :
\min(u,f_2) > z^{1/10} \}, \\
\N_5 &:=    \{ n\in \HH^* :
Z<\min(u,f_2) \le  z^{1/10} \}.
\end{align*}
In the above decomposition, the main parts are $\N_2$ and $\N_5$.  We
expect $\N_2$ to be small since, conditionally on $\bo(n;z)=k$, the
normal value of $\bo(n;z_h,z)$ is $hk/\log_2z>\ft{19}{10}h$.
It is more difficult to see that $\N_5$ is small too. This follows from the
fact that we count integers in this set
according to their number of factorizations in the form $n=uvf_1f_2$
with $y<vf_1<vf_2\le z$.  Suppose for instance that $f_1,f_2 \le z_j$.
For $\bo(n;z)=k$ and $\bo(n;z_j,z)=G$, then, ignoring the given 
information on the
localization of $vf_1$ and $vf_2$ in $(y,z]$, there are
$4^{k-G} 2^G = 4^k 2^{-G}$ such factorizations.
Thus, larger $G$
means fewer factorizations.  On probabilistic grounds, larger $G$
should also mean fewer factorizations when information on the 
localization of $vf_1$ and
$vf_2$ is available.
\par
\medskip
We now briefly consider the cases of $\N_0$ and $\N_1$.\par
Trivially,
\be\label{N0}
|\N_0| \le \frac{x}{\log z} + \sum_{d>\log z} \frac{x}{d^2} \ll
\frac{x}{\log z} \ll \frac{x}{(\log y)^{Q(\lam)} \sqrt{\log_2 y}},
\ee
since $Q(\lambda)\le Q(2)=\log 4-1$ in the range under consideration.
\par
By the argument on pages 40--41  of \cite{Divisors},
$$
\sum_{\substack{n\le x \\ \bo(n;z) > K}} 1 \ll
\frac{x}{(\log y)^{Q(\lam)} \sqrt{\log_2 y}}.
$$
Setting $t:=1-\ft32\eps$, Lemma \ref{Shiu} gives
\begin{align*}
\sum_{\substack{n\le x \\ \tau(n;y,z)\ge 1 \\ \bo(n;z) \le K_0
      }} 1 & \le t^{-(2-3\eps)\log_2 z} \sum_{\substack{dm\le x \\ y<d\le z}}
      t^{\bo(d) + \bo(m;z)} \ll x (\log z)^{2t-2-\b-(2-3\eps)\log t} \\
      &\ll x (\log y)^{-\b-2\eps^2} \ll x (\log
y)^{-Q(\lam)-\eps^2/2}.
\end{align*}
Therefore,
\be\label{N1}
|\N_1| \ll \frac{x}{(\log y)^{Q(\lam)} \sqrt{\log_2 y}}.
\ee

In the next four sections, we show that
\be\label{Nj}
|\N_j|\ll \frac{x}{(\log y)^{Q(\lam)}
\sqrt{\log_2 y}}\qquad (2\le j\le 5).
\ee
Together with \eqref{N0} and \eqref{N1}, this will complete the proof of
Theorem \ref{main_thm}.

%
\section{Estimation of $|\N_2|$}
%

We plainly have $|\N_2|\le \sum_{h}|\N_{2,h}|$.
For $1\le h\le 5\varepsilon\log_2z$, the numbers $n\in \N_{2,h}$ satisfy
$$\left\{\begin{array}{ll}
&x/\log z<n\le x,\\
\noalign{\vskip-3mm}\\
&k:=\bo(n;z)=K-b,\quad 0\le b\le 3\varepsilon\log_2z,\\
\noalign{\vskip-3mm}\\
&\bo(n;z_h,z)\le \ft{19}{10}h-\ft1{100}b,
\end{array}
\right.
$$
We note at the outset that $\N_{2,h}$ is empty unless
$h \ge b/190$.\par
Write $n=du$ with $y<d\le z$ and
$u\le x/y$. Let
$$
\bo(d;z_h) = D_1, \quad \Omega(d;z_h,z)=D_2, \quad
\bo(u;z_h)=U_1, \quad \Omega(u;z_h,z)=U_2,
$$
so that $D_1+D_2\ge 1$,
$D_2+U_2 \le \ft{19}{10}h - \ft1{100}b$ and $D_1+D_2+U_1+U_2=k$.

Fix $k=K-b$, $h$, $D_1$, $D_2$, $U_1$ and $U_2$.
By Lemma \ref{lem23} (with $w=z_h$, $a=U_1$, $b=U_2$),
the number of $u$ is
$$
\ll \frac{x}{y\log z} \; \frac{(\log_2 z-h)^{U_1}}{U_1!} (U_2+1)
\frac{(h+C)^{U_2}}{U_2!}.
$$
A second application of Lemma \ref{lem23} yields that the number of $d$ is
$$
\ll \frac{\eta y}{\log z} \frac{(\log_2 z - h)^{D_1}}{D_1!} (D_2+1)
\frac{(h+C)^{D_2}}{D_2!}.
$$
Since $D_2+U_2<2h$, we have $(h+C)^{U_2+D_2} \le \er^{2C} h^{U_2+D_2}$.
Summing over $D_1,D_2,U_1,U_2$ with $G=D_2+U_2$ fixed and using
the binomial theorem,
we find that the number of $n$ in question is
$$
\ll \frac{\eta x}{(\log z)^2} (\log_2 z-h)^{k-G} h^G (G+1)^2
\sum_{\substack{U_1+D_1=k-G
        \\ D_2+U_2 = G }} \frac{1}{U_1 ! D_1 ! D_2 ! U_2 !} \ll
\frac{\eta x 2^k}{(\log z)^2} A(h,G),
$$
where
$$
A(h,G) = (G+1)^2 \frac{(\log_2 z - h)^{k-G} h^G}{(k-G)! G!}.
$$
Since  $G+1\le G_h := \fl{\ft{19}{10}h}$, we have
$$
\frac{A(h,G+1)}{A(h,G)} \ge \frac{h(k-G)}{(G+1)(\log_2 z -h)} \ge
\frac{k-10\eps \log_2 z}{1.9 (1-5\eps)\log_2 z} > \ft{21}{20}
$$
if $\eps$ is small enough.  Next,
\begin{align*}
A(h,G_h) &\le (G_h+1)^2 \frac{(\log_2 z - h)^{k-G_h} (hk)^{G_h}}
      {k! (G_h/\er)^{G_h}} \\
&\ll (h+1)^2 \frac{(\log_2 z)^k}{k!}
      \big(\ft{20}{19}\er\big)^{19h/10} \er^{-h(k-G_h)/\log_2 z} \\
&\ll \frac{(\log_2 z)^k}{k!} \er^{-h/500},
\end{align*}
since $(k-G_h)/\log_2z>2-13\varepsilon$ and $\ft{19}{10}\log
(\ft{20}{19}\er)<2-1/400.$
Thus,
$$
\sum_{b/190 \le  h \le 5\eps \log_2 z} \;\; \sum_{0\le G \le G_h} A(h,G)
\ll \sum_{b/190 \le h \le 5\eps \log_2 z} A(h,G_h) \ll
\frac{(\log_2 z)^k}{k!} \er^{-b/95000}
$$
and so
$$\sum_{\di{n\in\N_2}{\bo(n;z)=k}}1\ll\frac{\eta x(2\log_2z)^k}{(\log
z)^2k!}\er^{-(K-k)/95000}\ll\frac{x\er^{-(K-k)/95000}}{(\log
y)^{Q(\lam)} \sqrt{\log_2 y}}, $$
by \eqref{factorial}. Summing over the range $K_0\le k\le K$
furnishes the required estimate \eqref{Nj}
for~$j=2.$

%
\section{Estimation of $|\N_3|$}
%

All integers $n=f_1f_2uv$ counted in $\N_3$ verify
$$\left\{
\openup10mm\begin{array}{ll}
&{x/\log z}<n\le x,\\
\bline\\
&\bo(n;z)\le K, \\
\bline\\
&y<vf_1<vf_2\le z,\quad \min(u,f_2)\le Z,\\
\end{array}
\right.
$$
where $Z$ is defined in \eqref{Z}. This is all we shall use in
bounding $|\N_3|$.
\par\goodbreak
Let $\N_{3,1}$ be the subset corresponding to the condition $f_2\le Z$ and let
$\N_{3,2}$ comprise those $n\in\N_3$ such that $u\le Z$.
\par
If $f_2 \le Z$, then $v>z^{1/2}$ and $u>x/\{vZ^2\log z\}>x^{1/3}$.
For $\frac12 \le t \le 1$ we have
\begin{align*}
|\N_{3,1}| &\le \sum_{f_1,f_2,v,u} t^{\bo(f_1 f_2 u v;z)-K} \\
&= t^{-K} \sum_{f_1 \le Z} t^{\bo(f_1)} \sum_{f_1<f_2\le \er^{\eta} f_1}
t^{\bo(f_2)} \sum_{y/f_1 < v \le z/f_1} t^{\bo(v)}
\sum_{u \le x/f_1f_2 v} t^{\bo(u;z)}.
\end{align*}
Apply Lemma \ref{Shiu} to the three innermost sums.  The $u$-sum  is
$$
\ll \frac{x}{f_1f_2 v} (\log z)^{t-1} \le \frac{x}{f_1 y} (\log z)^{t-1}.
$$
and the $v$-sum is
$$
\ll \frac{\eta y}{f_1} (\log z)^{t-1}.
$$
The  $f_2$-sum is $\ll \eta f_1 (\log f_1)^{t-1}$ if $f_1>\eta^{-3}$ and
otherwise is $\ll \eta f_1$ trivially (note that $\eta f_1\gg1$
follows from the fact that $(f_1+1)/f_1\le f_2/f_1\le
\er^\eta$).  Next
\begin{align*}
\sum_{f_1\le \eta^{-3}} \frac{1}{f_1} + \sum_{2\le f_1\le Z}
\frac{t^{\bo(f_1)}}{f_1} (\log  f_1)^{t-1} &\ll
\log_2 z + (\log_2 z) \max_{j\le \log_2 Z} \er^{j(t-1)}
\sum_{f_1\le \exp\{\er^j\}} \frac{t^{\bo(f_1)}}{f_1}  \\
&\ll (\log_2 z) (\log Z)^{2t-1}.
\end{align*}
Thus,
$$
|\N_{3,1}| \ll x (\log_2 x) (\log x)^E$$
with
$
E=-2\b - \lam \log t + 2t-2 + (2t-1)(1-4\eps).
$
We select optimally $t:=\ft14\lam/(1-2\eps)$, and check that $t
\ge\frac12$ since
$\lam \ge 2 -\eps/\log 2$.  Then
\begin{align*}
E &= -Q(\lam) + \lam \log(1-2\eps) + 4\eps \le -Q(\lam) + (2-\eps/\log
2)(-2\eps-2\eps^2) + 4\eps \\ &< -Q(\lam)-\eps^2.
\end{align*}

Next, we consider the case when $u\le Z$.  We observe that this implies
$$\ft14vz^2\le vx\le vn\log z=uf_1vf_2v\log z\le Zz^2\log z$$
hence $v\le 4Z\log z\le Z^2$, and therefore
$$\min(f_1,f_2)>z^{1/2}.$$  Also, $z>x^{1/3}$ since
$x/\log z<n=uvf_1f_2\le Zz^2.$
Thus, for $\frac12 \le t \le 1$, we have
\begin{align*}
|\N_{3,2}| &\le \sum_{f_1,f_2,v,u} t^{\bo(f_1 f_2 u v;z)-K} \\
&= t^{-K} \sum_{v\le Z^2} t^{\bo(v)} \sum_{u\le {xv/y^2}}
t^{\bo(u)}
\sum_{{y/v} < f_1 \le {z/v}} t^{\bo(f_1)}
\sum_{{y/v} < f_2 \le {z/v}} t^{\bo(f_2)}.
\end{align*}
The sums upon $f_1$ and $f_2$ are each
$$
\ll \frac{\eta y}{v} (\log z)^{t-1}
$$
and the $u$-sum is
$$
\ll \frac{xv}{y^2} (\log 2xv/y^2)^{t-1} \le \frac{xv}{y^2} (\log 2v)^{t-1}.
$$
Thus, selecting the same value $t:=\ft14\lam/(1-2\eps)$, we obtain
\begin{align*}
|N_{3,2}| &\ll t^{-K} x \eta^2 (\log z)^{2t-1} \sum_{v\le Z^2}
\frac{t^{\bo(v)} (\log 2v)^{t-1}}{v} \\
&\ll x(\log_2 z) (\log z)^E \le x (\log_2 z) (\log z)^{-Q(\lam)-\eps^2}.
\end{align*}
\par
This completes the proof of \eqref{Nj} with $j=3$.

%
\section{Estimation of $|\N_4|$}
%

We now consider those integers $n=f_1f_2uv$ such that
$$\left\{
\openup10mm\begin{array}{ll}
&{x/\log z}<n\le x,\\
\bline\\
&k:=\bo(n;z)=K-b,\quad 0\le b\le 3\varepsilon\log_2z,\\
\bline\\
&y<vf_1<vf_2\le z,\quad \min(u,f_2)>z^{1/10}.\\
\end{array}
\right.
$$
\par
With the notation \eqref{notnpf}, fix $k$, $F_1$, $F_2$, $U$ and $V$.
Here $u,f_1$ and $f_2$ are all
$>\frac12 z^{1/10}$.
By Lemma~\ref{lem23} (with $w=z$),
for each triple $f_1,f_2,v$ the number of $u$
is
$$
\ll \frac{x}{f_1 f_2 v \log z}\; \frac{(\log_2 z)^U}{U!}\cdot
$$
Using Lemma \ref{lem23} two more times, we obtain, for each $v$,
$$
\sum_{{y/v} < f_1 \le {z/v}} \frac{1}{f_1}
\sum_{{y/v} < f_2 \le {z/v}} \frac{1}{f_2}
\ll \frac{\eta^2}{(\log z)^2} \; \frac{(\log_2 z)^{F_1+F_2}}{F_1! F_2!}\cdot
$$
Now, Lemma \ref{lemA} gives
$$
\sum_v \frac{1}{v} \ll \frac{(\log_2 z)^V}{V!}\cdot
$$
Gathering these estimates and using \eqref{factorial} yields
\begin{align*}
|\N_4| & \ll \frac{x\eta^2}{(\log z)^3} \sum_{(2-3\eps)\log_2 z \le k \le K}
\;\; \sum_{F_1+F_2+U+V=k} \frac{(\log_2 z)^k}{F_1! F_2! U! V!} \\
&= \frac{x\eta^2}{(\log z)^3} \sum_{(2-3\eps)\log_2 z \le k \le K}
\frac{(2\log_2 z)^k}{k!} 2^k\\
&\ll \frac{x}{(\log y)^{Q(\lam)} \sqrt{\log_2 y}}\frac{2^K\eta}{\log z}\ll
\frac{x}{(\log y)^{Q(\lam)} \sqrt{\log_2 y}}\cdot
\end{align*}
Thus \eqref{Nj} holds for $j=4$.

%
\section{Estimation of $|\N_5|$}
%

It is plainly sufficient to bound the number of those  $n=f_1f_2uv$
satisfying the following conditions
$$\left\{
\openup10mm\begin{array}{ll}
&{x/\log z}<n\le x,
\\
\bline\\
&k:=\bo(n;z)=K-b,\quad 0\le b\le 3\varepsilon\log_2z,\\
\noalign{\vskip-3mm}\\
&\bo(n;z_h,z)> \ft{19}{10} h-\ft1{100}b\quad (1\le h\le 5\varepsilon\log_2z)\\
\noalign{\vskip-3mm}\\
&y<vf_1<vf_2\le z,\quad Z<\min(u,f_2)\le z^{1/10}.\\
\end{array}
\right.
$$

Define $j$ by $z_{j+2} < \min(u,f_2) \le z_{j+1}$. We have $1\le j\le
5\eps \log_2 z$.  Let $\N_{5,1}$ be the set of
those $n$ satisfying the above conditions with $u\le z_{j+1}$ and let
$\N_{5,2}$ be the complementary set, for which
$f_2\le z_{j+1}$.
\par
If $u\le z_{j+1}$, then
$v \le (z^2u\log z)/x \le 4u\log z \le z_j$ and $f_2>f_1>z^{1/2}$.
Recall notation \eqref{notnpf} and write
$$
F_{11}:=\bo(f_1;z_j), \quad  F_{12}:=\Omega(f_1;z_j,z), \quad
F_{21}:=\bo(f_2;z_j), \quad  F_{22}:=\Omega(f_2;z_j,z),
$$
so that the initial condition upon $\bo(n;z_h,z)$ with $h=j$ may be
rewritten as
$$F_{12}+F_{22} \ge G_j := \max(0,\fl{\ft{19}{10}j-{b/100}}).  $$
We count those $n$ in a dyadic interval
$(X,2X]$, where $x/(2\log z) \le X \le x$.
Fix $k,j,X,U,V,F_{rs}$ and apply Lemma \ref{lem23} to sums
over $u,f_1,f_2$.  The number of $n$ is question is
\begin{align*}
&\le \sum_{v\le z_j} \;\; \sum_{{vX/z^2} \le u \le {2vX/y^2}}
      \;\; \sum_{{y/v} < f_1 \le {z/v}}
      \;\; \sum_{{y/v} < f_2 \le {z/v}} 1 \\
&\ll \frac{\eta^2 X \er^j}{(\log z)^3} \frac{(\log_2 z -
      j)^{U+F_{11}+F_{21}}}{U! F_{11}! F_{21}!} (F_{12}+1)(F_{22}+1)
\frac{(j+C)^{F_{12}+F_{22}}}{F_{12}! F_{22}!}  \sum_{v\le z_j} \frac{1}{v}.
\end{align*}
Bounding the $v$-sum by Lemma \ref{lemA}, and summing over $X$,
$U,V,F_{rs}$ with $F_{12}+F_{22}=G$ yields
$$
|\N_{5,1}| \ll \frac{\eta^2 x}{(\log z)^3} \sum_{(2-3\eps)\log_2 z \le k \le
      K} 4^k \sum_{1\le j\le 5\eps\log_2 z} \;\; \sum_{G_j \le G \le k} M(j,G),
$$
where
$$
M(j,G):=\er^j (G+1)^2 \frac{(\log_2 z - j)^{k-G} (j+C)^G}{2^G (k-G)! G!}.
$$
Let $j_b=\fl{\ft12b + 100C+100}$.  If $j\le j_b$, then
$j+C \le \ft{99}{100} (j+C_b)$ with $C_b:=3C+2+\frac{b}{100}$ and,
introducing $R:=\max_{G\ge 0} \big\{(G+1)^2
(\ft{99}{100})^G\big\}$, we have
\begin{align*}
\sum_{1\le j\le j_b} \;\; \sum_{G_j\le G\le k} M(j,G)&\le R\sum_{1\le j\le j_b}
\er^j  \dsp\sum_{0\le G\le k} \frac{(\log_2 z-j)^{k-G} (j+C_b)^G}{2^G
      G! (k-G)!} \\
&\ll\frac{1}{k!} \sum_{1\le j\le j_b} \er^j \( \log_2 z - \ft12j
+\ft12 C_b \)^k \\
&\ll \frac{(\log_2 z)^k}{k!} \sum_{1\le j\le j_b}
\er^{j+({b/200}-{j/2}) {k/\log_2 z}} \\
&\ll \frac{(\log_2 z)^k}{k!} \er^{{b/100} + 2\eps j_b} \ll
\frac{(\log_2 z)^k}{k!} \er^{{b/50}}.
\end{align*}
When $j> j_b$, then $$G_j\ge
\ft95(j+C)+\ft1{10}(j_b+C+1)-\ft{1}{100}b-1\ge \ft95(j+C)+9\ge 189.$$
Thus, for
$G\ge G_j$ we have
$$
\frac{M(j,G+1)}{M(j,G)} = \pfrac{G+2}{G+1}^2 \frac{j+C}{2(G+1)}
\frac{k-G}{\log_2 z - j} \le \ft{4}{7}.
$$
Therefore,
\begin{align*}
\sum_{G_j \le G \le k} M(j,G) &\ll M(j,G_j) \ll \frac{j^2 \er^j}{k!}
      \frac{(\log_2 z-j)^{k-G_j} (jk)^{G_j}}{2^{G_j} G_j!} \\
&\le \frac{j^2\er^{j} (\log_2 z)^k}{k!} \er^{ -j(k-G_j)/\log_2 z}
\pfrac{\er jk}{2G_j \log_2 z}^{G_j}
\ll \frac{(\log_2 z)^k}{k!} \er^{-{j/5}},
\end{align*}
since $k-G_j\ge (2-10\eps)\log_2z$, $\er jk/(2G_j\log_2z)\le
\ft{5}{9}\er$, and $-1+\ft{19}{10}\log
(\ft{5}{9}\er)<-\ft1{5}$. We conclude that
\be\label{MhG}
\sum_{1\le j \le 5\eps \log_2 z} \;\; \sum_{G_j \le G \le k} M(j,G) \ll
\frac{(\log_2 z)^k}{k!} \er^{b/50}
\ee
and hence, by \eqref{factorial},
$$
|\N_{5,1}| \ll \frac{\eta^2 x}{(\log z)^3} \sum_{k\le K} \frac{(2\log_2
      z)^k}{k!} 2^{K-b/2} \ll \frac{\eta^2 2^K x}{(\log z)^3} \frac{(2\log_2
      z)^{K}}{K!} \ll \frac{x}{(\log y)^{Q(\lam)} \sqrt{\log_2 y}}.
$$

Now assume $f_2 \le z_{j+1}$.  Then $\min(u,v)>\sqrt{z}$.
Fix $F_1$, $F_2$ and
$$
\bo(v;z_j)=V_1, \quad  \Omega(v;z_j,z)=V_2, \quad  \bo(u;z_j)=U_1, \quad
\Omega(u;z_j,z)=U_2.
$$
By Lemma \ref{lem23}, given $f_1,f_2$ and $v$, the number of $u$ is
$$
\ll \frac{x}{f_1 f_2 v \log z} \frac{(\log_2 z-j)^{U_1} (U_2+1)
      (j+C)^{U_2}}{U_1! U_2!}.
$$
Applying Lemma \ref{lem23} again, for each $f_1$ we have
$$
\sum_{\di{f_1 < f_2 \le \er^\eta f_1}{{y/f_1} < v \le
        {z/f_1}}} \frac{1}{f_2 v} \ll \frac{\eta^2 \er^j}{(\log z)^2}
\frac{(V_2+1) (\log_2 z-j)^{V_1+F_2}(j+C)^{V_2}}{V_1! V_2! F_2!}.
$$
By Lemma \ref{lemA},
$$
\sum_{f_1\le z_j} \frac{1}{f_1} \ll \frac{(\log_2 z-j)^{F_1}}{F_1!}.
$$
Combine these estimates, and sum over $F_1,F_2,U_1,U_2,V_1,V_2$ with
$V_2+U_2=G$.  As in the estimation of $|\N_{5,1}|$, sum over
$k,j,G$ using \eqref{factorial} and \eqref{MhG}.  We obtain
\begin{align*}
|\N_{5,2}| &\ll \frac{\eta^2 x}{(\log z)^3} \sum_{(2-3\eps)\log_2 z \le k\le
K} 4^k \sum_{1\le j\le 5\eps\log_2 z} \;\;\sum_{G_j <G\le k} M(j,G) \\
&\ll \frac{x}{(\log y)^{Q(\lam)} \sqrt{\log_2 y}}.
\end{align*}

%
%

\bibliographystyle{amsplain}
\bibliography{ft9}

\end{document}